\documentclass[12pt]{article}
\usepackage{amsmath}
\usepackage{amssymb}
\usepackage{amsfonts}
\usepackage{amscd}
\usepackage{enumerate}
\usepackage{fancyheadings}
\usepackage{enumerate}
\usepackage{color}
\usepackage{graphicx}
\usepackage{amsthm}
\newtheorem{Lem}{Lemma}
\newtheorem{Th}{Theorem}

\newtheorem{Def}{Definition}

\newtheorem{Assume}{Assumption}

\DeclareMathOperator{\supp}{supp}

\DeclareMathOperator{\graph}{graph}

\DeclareMathOperator{\dist}{dist}
\DeclareMathOperator{\sgn}{sgn}

\begin{document}
\begin{center}
\bf{\Large Semi-Classical Behavior of the Spectral Function}\footnote{Prepared using AMS-LaTeX.}
\end{center}

\noindent
{\bf Ivana Alexandrova}

\noindent
Department of Mathematics, University of Toronto, Toronto, Ontario, Canada 
M5S 3G3, Tel.: 1-416-946-0318, Fax: 1-416-978-4107, email: alexandr@math.toronto.edu

\noindent
February 25, 2005

\begin{abstract}
We study the semi-classical behavior of the spectral function of the Schr\"{o}dinger operator with short range potential.
We prove that the spectral function is a semi-classical Fourier integral operator quantizing the forward and backward Hamiltonian flow relations of the system.
Under a certain geometric condition we explicitly compute the phase in an oscillatory integral representation of the spectral function.
\end{abstract}

\section{Introduction}
We study the structure of the spectral function associated with the
semi-classical Schr\"{o}dinger operator with short range potential on $\mathbb{R}^{n}.$
We prove that the appropriately cut-off spectral function is a semi-classical Fourier 
integral operator associated to the union of the backward and the forward 
Hamiltonian flow 
relations of the principle symbol of the 
operator.
We also show how this allows us, under a certain geometric assumption, to 
compute the phase in an oscillatory integral representation of the spectral function.

Our result is motivated by the following theorem by Vainberg.
In \cite[Theorem XII.5]{V} Vainberg considers operators of the form
\begin{equation*}
A=\sum_{i, j=1}^{n}a_{i,j}\frac{\partial^{2}}{\partial x_i \partial 
x_j}+\sum_{i=1}^{n}b_i\frac{\partial}{\partial x_i}+c,
\end{equation*}
where $a_{i,j}, b_{i}, c\in C^{\infty}(\mathbb{R}^{n})$ and $A\equiv -\Delta$ for $\|x\|\geq 
1.$ 
He assumes that the energy $1$ is non-trapping for the principal symbol of 
$A.$
Vainberg then establishes an asymptotic expansion in $\lambda\to\infty$ 
of the spectral function $e_{\lambda},$ which is defined as the 
Schwartz kernel of $\frac{d E_{\lambda}}{d\lambda},$ where 
$\{E_{\lambda}\}$ denotes the spectral 
family of $A.$ 
Vainberg expresses this asymptotic expansion in the form of a Maslov canonical 
operator $K_{\Lambda, \lambda}$ associated to a certain Lagrangian 
submanifold $\Lambda=\Lambda_y\subset T^{*}\mathbb{R}^{n}$ and acting on 
another asymptotic sum in $\lambda.$
The Lagrangian submanifold $\Lambda_y$ consists of 
the phase trajectories at energy 1 of the principal symbol of $A$ passing 
through a fixed base point 
$x(0)=y,$ while the terms of the asymptotic sum on which $K_{\Lambda, 
\lambda}$ acts solve a recurrent system of transport equations along the 
phase trajectories of the system.

\section{The Semi-Classical Spectral Function}
Here we study the semi-classical behavior of the spectral function of a 
Schr\"{o}dinger operator with short range potential at a 
fixed energy $\lambda>0.$
More precisely, we work in the following setting.
Let $X$ be a smooth manifold of dimension 
$n>1$ such that
$X$ coincides with $\mathbb{R}^{n}$ on 
$\mathbb{R}^{n}\backslash B(0, R_0)$ for some 
$R_{0}>0,$ where $B(0, R_0)=\{x\in\mathbb{R}^{n}: \|x\|<R_0\}$ and 
$\|\cdot\|$ denotes the Euclidean norm on $\mathbb{R}^{n}.$ 
Let $g$ be a Riemannian metric on $X$ which satisfies 
the condition \begin{equation*}
g_{ij}\left(x\right)=\delta_{ij} \text{ for } \|x\|>R_{0}.
\end{equation*}
Let $V\in C^{\infty}(X; \mathbb{R})$ be such that
\begin{equation*}
\left|\frac{\partial^{\alpha}}{\partial 
x^{\alpha}}V(x)\Big|_{\mathbb{R}^{n}\backslash B(0, R_0)}\right|\leq C_{\alpha}(1+\|x\|)^
{-\mu-|\alpha|},\; x\in\mathbb{R}^{n}, \mu>0.
\end{equation*}
Then the operators $P\left(h\right)=\frac{1}{2}h^{2}\Delta_{g}+V,$ 
$0<h\leq 1,$ admit unique self-adjoint extensions with common domain $H^{2}(X).$
We denote by $\{E_{\mu}(h)\}$ the spectral family of the operator $P(h)$ 
and by $e_{\mu}(h)$ the spectral function of $P(h),$ i.e. the Schwartz kernel of $\frac{dE_{\mu}(h)}{d\mu}.$
We also let 
$p(x, \xi)=\frac{1}{2}\|\xi\|^{2}_{g}+V(x)$ be the semi-classical 
principal symbol of $P.$

We set $R(\lambda\pm i0, h)=\lim_{\epsilon\downarrow 0}\left(P(h)-\lambda\mp 
i\epsilon\right)^{-1},$ where 
the 
limit is taken in the space of bounded operators $\mathcal{B}(L^{2}_{\alpha}(X), 
L^{2}_{-\alpha}(X)),$ $\alpha>\frac{1}{2},$ where 
\[L^{2}_{\pm\alpha}(X)=\left\{f\in L^{2}(X): 
f|_{B(0, R_0)}\langle\cdot\rangle^{\pm\alpha}\in 
L^{2}(\mathbb{R}^{n}\backslash B(0, R_0))\right\}.\]
The operator norm in $\mathcal{B}(L^{2}_{\alpha}(X), L^{2}_{-\alpha}(X))$ will be denoted by 
$\|\cdot\|_{\alpha, -\alpha}.$

Let $H_{p}$ be the Hamiltonian vector field of $p$ and let $\gamma(\cdot; x_0,
\xi_0)=(x\left(\cdot; x_0, \xi_0\right), \xi(\cdot; x_0, \xi_0))$ denote the integral curve
of $H_{p}$, or (phase) trajectory, with initial conditions $(x_0, \xi_0)\in T^{*}X.$
We define a non-trapping energy level as follows:
\begin{Def}
The energy $\lambda>0$ is non-trapping if for every $\left(x_{0}, 
\xi_{0}\right)$ with $\frac12\|\xi_0\|^{2}+V(x_0)=\lambda$ there exists $t_0>0$ such that 
$x\left(s; x_{0}, \xi_{0}\right)\in
\mathbb{R}^{n}\backslash B(0, R_0)$ for every $|s|>t_0.$
A phase trajectory $\gamma(\cdot; x_0, \xi_0)$ is non-trapped if there exists $t>0$
such that $x(s; x_{0}, \xi_{0})\in \mathbb{R}^{n}\backslash B(0, R_{0})$
for all $|s|>t.$
\end{Def}

We now choose functions $\chi_j\in C_{c}^{\infty}(\mathbb{R}^{n}\backslash B(0, R_0);
\mathbb{R}),$ $j=1, 2,$ with disjoint supports.
We assume that $\lambda>0$ is such that $P-\lambda$ is of principal type.
Then it follows that $\Sigma_{\lambda}=p^{-1}(\lambda)$ is a $2n-1-$ dimensional submanifold 
of 
$T^{*}\mathbb{R}^{n}$ and 
\[\Lambda^{+}_{R}(\lambda)=\left\{(y, -\eta; x, \xi): (y, 
\eta)\in\Sigma_{\lambda}, (x, \xi)=\exp(tH_p)(y, \eta), t>0\right\}\cap 
T^{*}(\supp\chi_1\times\supp\chi_2)\]
and
\[\Lambda^{-}_{R}(\lambda)=\left\{(x, -\xi; y, \eta): (x, 
\xi)\in\Sigma_{\lambda}, (y, \eta)=\exp(tH_p)(x, \xi), t<0\right\}\cap 
T^{*}(\supp\chi_1\times\supp\chi_2)\]
are Lagrangian submanifolds of $T^{*}\mathbb{R}^{n}\times T^{*}\mathbb{R}^{n}.$

To state our main theorem, we further let $\pi_2: T^{*}\mathbb{R}^{n}\times 
T^{*}\mathbb{R}^{n}\to T^{*}\mathbb{R}^{n}$ denote the canonical projection onto the second factor.
We also refer the reader to the Appendix for the definition of the class of semi-classical 
Fourier integral operators $\mathcal{I}_{h}^{r},$ $r\in\mathbb{R},$ as well as for a 
review of the relevant notions from semi-classical analysis.

We can now state the following
\begin{Th}\label{tspffio}
Let $\rho_0\in\Lambda_{R}^{+}(\lambda)$ be such that 
$\gamma\left(\pi_2\left(\rho_0\right)\right)$ is 
non-trapped.
Let $\|R(\lambda\pm i0, h)\|_{\alpha, -\alpha}=\mathcal{O}(h^{s}),$ $s\in\mathbb{R}.$

Then there exist open sets $W_{\pm}\in T^{*}\mathbb{R}^{n}\times 
T^{*}\mathbb{R}^{n},$ such that 
\[\chi_2 \frac{d E_{\lambda}}{d \lambda}\chi_1\in\mathcal{I}_{h}^{1}\left(\mathbb{R}^{2n}, 
\left(\overline{W}_{+}\cap\Lambda_{R}^{+}(\lambda)\right)\cup\left(\overline{W}_{-}\cap\Lambda_{R}^{-}(\lambda)\right)\right).\]
\end{Th}

Before proving the theorem, we would like to make two remarks:

\noindent
{\bf Remark 1.} The assumption on the polynomial bound of the resolvent is satisfied in a 
number of interesting situations: at non-trapping energies (see \cite[Lemma 2.2]{RT}) and at 
trapping energies $\lambda$ when we assume that the resonances 
$\left(\lambda_j\right)$ of 
$P(h)$ are such that $\left|\Im\lambda_j\right|\geq Ch^{q},$ if 
$\Re\lambda_j\in\left[\lambda-\epsilon, \lambda+\epsilon\right]$ for some $\epsilon>0$ (see \cite[Proposition 5.1]{M}).
In the latter case, in order to define the resonances by complex scaling, 
the author also assumes that there exists 
$\theta_0\in[0, \pi),$ $\epsilon>0,$ and $R>0$ such that $V$ extends holomorphically to 
\[D_{\epsilon, R, \theta_0}=\left\{r\omega: \omega\in\mathbb{C}^{n}, 
\dist(\omega, 
\mathbb{S}^{n-1})<\epsilon, 
r\in\mathbb{C}, |r|>R, \arg r\in[-\epsilon, \theta_0+\epsilon)\right\}\]
and
\[\exists\beta>0, \exists M>0, \forall x\in D_{\epsilon, R, \theta_0}, 
|V(x)|\leq 
C|x|^{-\beta}.\]  

\noindent
{\bf Remark 2.}  \cite[Theorem 1]{AIfio} roughly says that semi-classical Fourier integral 
distributions, i.e., Schwartz kernels of the elements of 
$\mathcal{I}_{h}^{r}(\mathbb{R}^{k}, \Lambda),$ $r\in\mathbb{R},$ $k\in\mathbb{N},$ where $\Lambda\subset 
T^{*}\mathbb{R}^{k}$ is a smooth closed Lagrangian submanifold, can be represented microlocally as oscillatory integrals $\int 
e^{i\phi(x, \theta)/h} a_{\phi}(x, \theta)d\theta$ with $a_{\phi}\in S(1)$ for every 
non-degenerate phase function $\phi=\phi(x, \theta)$ which parameterizes $\Lambda$ in the 
sense that 
$\Lambda=\{(x, d_x \phi): d_\theta \phi=0\}$ near some point $\rho\in\Lambda.$
Furthermore, such a phase function always exists near any point 
$\rho\in\Lambda$ (see \cite[Section 4.1]{AIfio}) and the corresponding symbol $a_{\phi}$ 
with an asymptotic expansion in $h$ can always be found (see \cite[Theorem 1]{AIfio}). 
Theorem \ref{tspffio} thus implies that the appropriately cut-off spectral function always 
admits 
such an oscillatory integral representation.
In Lemma \ref{laction} below we explicitly compute the phase for certain Lagrangians.

\begin{proof}[Proof of Theorem \ref{tspffio}]
We recall that 
\begin{equation*} 
\chi_2 R(\lambda+i0, h)\chi_1-\chi_2 R(\lambda-i0, h)\chi_1=2\pi 
i\chi_2\frac{d E_{\lambda}}{d \lambda}\chi_1
\end{equation*}

Since $\gamma(\pi_{2}(\rho_0))$ is non-trapped, it follows that there
exists
an open set $W_{+}\subset T^{*}\mathbb{R}^{n}\times T^{*}\mathbb{R}^{n}$ 
such that $\gamma(\pi_{2}(\rho))$ is non-trapped for every $\rho \in \overline{W}_{+}.$
Let $W_{-}=\{(x, -\xi; y, \eta): (y, -\eta; x, \xi)\in W_{+}\}.$

By \cite[Lemma 1]{AI} and the estimate $\|R(\lambda\pm i0, h)\|_{\alpha, 
-\alpha}=\mathcal{O}(h^{s})$ we then obtain that $K_{\chi_2 R(\lambda\pm i0, h)\chi_1}\in\mathcal{D}'_{h}(\mathbb{R}^{2n}),$ 
where $K_{\chi_2 R(\lambda\pm i0, h)\chi_1}$ denotes the Schwartz kernel of $\chi_2 
R(\lambda\pm i0, h)\chi_1.$
The same proof as in \cite[Theorem 2]{AI} then shows that 
\[\chi_2 R(\lambda\pm i0, h)\chi_1\in\mathcal{I}_{h}^{1}\left(\mathbb{R}^{2n}, 
\overline{W}_{\pm}\cap\Lambda^{\pm}_{R}(\lambda)\right).\] 
The assumption that $\chi_1$ and $\chi_2$ have disjoint supports is 
essential in the proof of \cite[Theorem 2]{AI}.
Since $\Lambda^{+}_{R}(\lambda)$ and $\Lambda^{-}_{R}(\lambda)$ are disjoint, it follows 
that 
\begin{equation*}
\chi_2 \frac{d E_{\lambda}}{d\lambda}\chi_1\in\mathcal{I}_{h}^{1}\left(\mathbb{R}^{2n},
\left(\overline{W}_{+}\cap\Lambda^{+}_{R}(\lambda)\right)\cup
\left(\overline{W}_{-}\cap\Lambda^{-}_{R}(\lambda)\right)\right).
\qedhere
\end{equation*}
\end{proof}

We now turn to showing how the forward and backward flow relations, 
$\overline{W}_{+}\cap\Lambda^{+}_{R}(\lambda)$ and 
$\overline{W}_{-}\cap\Lambda^{-}_{R}(\lambda),$ respectively, can be 
parameterized 
by a non-degenerate phase function.
For that we make the following assumption
\begin{Assume}\label{cf}
The trajectory $\gamma_0(\cdot; y_0, \eta_0)=(x_0(\cdot; y_0, \eta_0), 
\xi_0(\cdot; y_0, \eta_0))\subset\Sigma_{\lambda}$ with 
$y_0\in\supp\chi_1$ and $x_0(t_0; y_0, \eta_0)=z_0\in\supp\chi_2$ is 
non-trapped and is
contained in a central field, i.e. (see \cite[Section 46.C]{Arnold}),
\[\det\left(\frac{\partial x_0}{\partial \eta}\left(t_0; y_0, 
\cdot\right)\right)(\eta_0)\ne 0.\]
\end{Assume}
By the Implicit Function Theorem, then, there exist open neighborhoods 
$T\times Y\times Z$
of
$t_0, y_0, x_0$ and a unique function $\eta\in C^{\infty}(T\times Y\times Z)$ such
that $\eta(t_0, y_0, x_0)=\eta_0$ and $x(t; y, \eta(t, y, x))=x.$
We can therefore define the action 
\begin{equation*}
S_{\nu}(y, z, t)=\int_{l(t, y, z)}L(\dot{x}, x)dt
\end{equation*}
over the segment $l(t, y, z)$ from $(y, \eta(t, y, z))$ to $(z, \xi(t; y, 
\eta(t, y, z)))$ of the trajectory $\gamma(y, \eta(t, y, z))$ for 
$(t, y, z)\in T\times Y\times Z,$ where $L(\dot{x}, 
x)=\frac12\|\dot{x}\|^{2}-V(x)+\lambda$ is the Lagrangian 
associated to the Hamiltonian $p-\lambda$ and $\nu=\sgn t_0.$ 

We now have the following

\begin{Lem}\label{laction}
Let $\lambda>0$ be such that $P-\lambda$ is of principal type and let $\gamma_0,$ $T,$ $Y,$ 
$Z,$ $\nu$ be as above.

Then $S_\nu$ is a non-degenerate phase function and 
$\Lambda_{S_{\nu}}=\{(z, y, d_z S_\nu, d_y S_\nu):d_t S_\nu=0, (t, y, 
z)\in T\times Y\times 
Z\}$ 
is a closed Lagrangian 
submanifold such that 
$\Lambda_{S_{\nu}}\subset\Lambda_{R}^{\nu}(\lambda).$ 
\end{Lem}

\begin{proof}
Assumption \ref{cf} allows us to apply \cite[Theorem 46.C]{Arnold} and we obtain 
\begin{equation}\label{dz}
d_z S_{\nu}(t; y, z)=\xi(t; y, \eta(t, y, z)).
\end{equation}
From
\[\det\left(\frac{\partial x}{\partial \eta}\left(t_0; y_0, \cdot\right)\right)(\eta_0)\ne 
0\]
it follows that there exists an open set $U\subset T^{*}\mathbb{R}^{n}$ with $(y_0, 
\eta_0)\in U$
 such that
\[\pi: \{(y, \eta; \exp(t_0 H_p)(y, \eta)): (y, \eta)\in U\}\rightarrow \mathbb{R}^{n}\times
\mathbb{R}^{n},\]
where $\pi:T^{*}\mathbb{R}^{n}\times T^{*}\mathbb{R}^{n}\rightarrow \mathbb{R}^{n}\times
\mathbb{R}^{n}$ is the canonical projection, is a diffeomorphism.
Therefore there exists $\psi\in C^{\infty}(\mathbb{R}^{2n})$ such that
$\graph\exp(t_0 H_p)=\Lambda_{\psi}$ near $(y_0, \eta_0; z_0, \xi(t_0; y_0, \eta_0)).$
Since $\exp(t_0 H_p)$ is a symplectomorphism, it follows that
$\det\left(\frac{\partial^{2}\psi(\cdot, \cdot\cdot)}{\partial_x \partial_y}\right)
(x_0, y_0)\ne 0.$
This implies that $\det\left(\frac{\partial \xi}{\partial y}(t_0; \cdot, \eta_0)\right)(y_0)
\ne 0.$
This, on the other hand, implies that
\begin{equation*}
\det\left(\frac{\partial x}{\partial \eta}(-t_0; x_0,
\cdot)\right)(\xi_0)\ne 0, \text{ where } (x_0, \xi_0)=\exp(t_0 H_p)(y_0, \eta_0).
\end{equation*}
Thus we can again apply \cite[Theorem 46.C]{Arnold} and obtain
\begin{equation}\label{dy}
d_y S_{\nu}(t; y, z)=-\eta(t, y, z).
\end{equation}

Lastly, we observe that since $P-\lambda$ is of principal type, it follows
that $dd_t S_{\nu}\ne 0$ on $\{d_t S_{\nu}=p-\lambda=0\}$ and therefore
$S_{\nu}$ is a
non-degenerate phase function.
This, together with \eqref{dz} and \eqref{dy}, implies that, perhaps after decreasing 
$T\times Y\times Z$ 
around $(t_0, y_0, 
z_0),$
\begin{equation*} 
\begin{aligned}
\Lambda_{S_{\nu}}=\bigg\{&\left(z, y, \xi\left(t; y, \eta\left(t, y, z\right)\right), \eta\left(t, y, z\right)\right):\\ 
& d_t 
S_\nu(t, y, z)=\frac{1}{2}\left\|\dot{x}\left(t; y, \eta\left(t, y, z\right)\right)\right\|^{2}+V\left(x\left(t; y, \eta\left(t, y, z\right)\right)\right)-\lambda=0,\\
& (t, y, z)\in T\times Y\times X\bigg\}\subset\Lambda_{R}^{\nu}(\lambda)
\end{aligned}
\end{equation*}
and $\Lambda_{S_{\nu}}$ is a closed Lagrangian submanifold of $T^{*}\mathbb{R}^{2n}.$
\end{proof}

The following lemma describes the microlocal structure of the cut-off spectral function.

\begin{Lem}
Let $\lambda>0$ be such that $P-\lambda$ is of principal type and $\|R(\lambda+i0, 
h)\|_{\alpha, -\alpha}=\mathcal{O}(h^s),$ $s\in\mathbb{R}.$
Let $\gamma(\cdot; y_0, \eta_0)$ with $(y_0, \eta_0)\in\Sigma_\lambda,$ 
$y_0\in\supp\chi_1,$ and $x(t_0; y_0, \eta_0)=z_0\in\supp\chi_2$ be 
non-trapped and contained in a central field. 
Let $\nu=\sgn t_0.$

Then there exists $a_{\nu}\in S_{2n+1}^{\frac{n+3}{2}}(1)\cap 
C_{c}^{\infty}(\mathbb{R}^{2n+1})$ such that \[\left(\chi_2\otimes\chi_1\right)e_\lambda=\int 
e^{\frac{i}{h}S_{\nu}(y, x, t)}a_{\nu}(y, x, t)dt\] microlocally near 
$(y_0, -\eta_0; z_0, \xi(t_0; y_0, \eta_0))$ if $t_0>0$ or near $(z_0, 
-\xi(t_0; y_0, \eta_0); y_0, \eta_0)$ if $t_0<0.$
\end{Lem}

\begin{proof}
With Theorem \ref{tspffio} and Lemma \ref{laction} the conditions of 
\cite[Theorem 1]{AIfio} are satisfied and we obtain that there exist 
$a_{\nu}\in S_{2n+1}^{\frac{n+3}{2}}(1)\cap 
C_{c}^{\infty}(\mathbb{R}^{2n+1})$ such that 
\[\left(\chi_2\otimes\chi_1\right)e_\lambda=\int e^{\frac{i}{h}S_{\nu}(y, 
x, t)}a_{\nu}(y, x, t)dt\] microlocally near 
$(y_0, -\eta_0; z_0, \xi(t_0; y_0, \eta_0))$ for $t_0>0$ and near $(z_0, 
-\xi(t_0; y_0, \eta_0); y_0, \eta_0)$ for $t_0<0.$
\end{proof}

\appendix
\section{Elements of Semi-Classical Analysis}
In this section we recall some of the elements of semi-classical analysis
which we
will use here.
First we define two classes of symbols
\begin{equation*}
S_{2n}^{m}\left(1\right)= \left\{ a\in
C^{\infty}\left(\mathbb{R}^{2n}\times(0, h_0]\right): \forall
\alpha, \beta\in\mathbb{N}^{n}, \sup_{(x, \xi,
h)\in\mathbb{R}^{2n}\times (0,
h_{0}]}h^{m}\left|\partial^{\alpha}_{x}\partial^{\beta}_{\xi}a\left(x,
\xi;
h\right)\right|\leq
C_{\alpha, \beta}\right\}
\end{equation*}
and
\begin{equation*}
S^{m, k}\left(T^{*}\mathbb{R}^{n}\right)=\left\{a\in
C^{\infty}\left(T^{*}\mathbb{R}^{n}\times(0, h_0]\right): \forall \alpha,
\beta\in\mathbb{N}^{n}, \left|\partial^{\alpha}_{x}\partial^{\beta}_{\xi}
a\left(x, \xi;
h\right)\right|\leq
C_{\alpha,
\beta}h^{-m}\left\langle\xi\right\rangle^{k-|\beta|}\right\},
\end{equation*}
where $h_0\in(0,1]$ and $m, k\in\mathbb{R}.$
For $a\in S_{2n}\left(1\right)$ or $a\in S^{m, 
k}\left(T^{*}\mathbb{R}^{n}\right)$ we define
the
corresponding semi-classical pseudodifferential operator of class
$\Psi_{h}^{m}(1, \mathbb{R}^{n})$ or $\Psi_{h}^{m, k}(\mathbb{R}^{n}),$ 
respectively, by
setting
\begin{equation*}
Op_{h}\left(a\right)u\left(x\right)=\frac{1}{\left(2\pi
h\right)^{n}}\int\int e^{\frac{i\left\langle x-y,
\xi\right\rangle}{h}}a\left(x, \xi; h\right)u\left(y\right) dy d\xi, 
\;u\in
\mathcal{S}\left(\mathbb{R}^{n}\right),
\end{equation*}
and extending the definition to $\mathcal{S}'\left(\mathbb{R}^{n}\right)$ 
by
duality.
Below we shall work only with symbols which admit asymptotic expansions in 
$h$ and with
pseudodifferential operators which are quantizations of such symbols.
For $A\in\Psi_{h}^{k}(1, \mathbb{R}^{n})$ or $A\in\Psi_{h}^{m, 
k}(\mathbb{R}^{n}),$ we shall
use $\sigma_{0}(A)$ and $\sigma(A)$ to denote its principal symbol and its 
complete symbol,
respectively.
A semi-classical pseudodifferential operator will be called of principal
type if its
principal symbol $a_0$ satisfies
\begin{equation*}
a_0=0\implies da_0\ne 0.
\end{equation*}

We also define the class of semi-classical distributions
$\mathcal{D}_{h}'(\mathbb{R}^{n})$ with which we will work here
\begin{equation*}
\begin{aligned}
\mathcal{D}'_{h}(\mathbb{R}^{n}) = & \big\{u\in C^{\infty}_{h}\left((0,
1];
\mathcal{D}'\left(\mathbb{R}^{n}\right)\right): \forall\chi\in
C_{c}^{\infty}\left(\mathbb{R}^{n}\right) \exists\: N\in\mathbb{N}\text{
and
} C_{N}>0:\\
& \quad |\mathcal{F}_{h}\left(\chi u\right)\left(\xi, h\right)|\leq
C_{N}h^{-N}\langle\xi\rangle^{N}\big\}
\end{aligned}
\end{equation*}
where
\begin{equation*}
\mathcal{F}_{h}\left(u\right)\left(\xi,
h\right)=\int_{\mathbb{R}^{n}}e^{-\frac{i}{h}\left\langle x,
\xi\right\rangle}u\left(x, h\right)dx
\end{equation*}
with the obvious extension of this definition to
$\mathcal{E}_{h}'(\mathbb{R}^{n}).$
Everywhere here we work with the $L^{2}-$ based semi-classical Sobolev spaces 
$H^{s}(\mathbb{R}^{n}),$ $s\in\mathbb{R},$ which consist of the distributions 
$u\in\mathcal{D}_{h}'(\mathbb{R}^{n})$ such 
that 
\[\|u\|_{H^{s}(\mathbb{R}^{n})}^{2}=\frac{1}{(2\pi 
h)^{n}}\int_{\mathbb{R}^{n}}(1+\|\xi\|^{2})^{s}\left|\mathcal{F}_{h}(u)(\xi, 
h)\right|^{2}d\xi<\infty.\]

We shall say that $u=v$ {\it microlocally} (or $u\equiv v$) near an open 
or closed set
$U\subset T^{*}\mathbb{R}^{n}$, if
$P(u-v)=\mathcal{O}\left(h^{\infty}\right)$ in
$C_{c}^{\infty}\left(\mathbb{R}^{n}\right)$ for
every $P\in \Psi^{0}_{h}\left(1, \mathbb{R}^{n}\right)$ such that
\begin{equation*}
WF_{h}\left(P\right)\subset \tilde{U}, \bar{U}\Subset \tilde{U}\Subset
T^{*}\mathbb{R}^{n}, \tilde{U} \text{ open}.
\end{equation*}
We shall also say that $u$ satisfies a property $\mathcal{P}$  {\it
microlocally} near an open set $U\subset T^{*}{\mathbb{R}^{n}}$ if there
exists $v\in\mathcal{D}_{h}'\left(\mathbb{R}^{n}\right)$ such that $u=v$
microlocally
near $U$ and $v$ satisfies property $\mathcal{P}$.

For open sets $U, V\subset T^{*}\mathbb{R}^{n},$ the operators $T,
T'\in\Psi^{m}_{h}\left(\mathbb{R}^{n}\right)$ are said to be {\it
microlocally
equivalent} near $V\times U$ if for any $A,
B\in\Psi_{h}^{0}\left(\mathbb{R}^{n}\right)$
such that
\begin{equation*}
WF_{h}\left(A\right)\subset\tilde{V},
WF_{h}\left(B\right)\subset\tilde{U},
\bar{V}\Subset\tilde{V}\Subset T^{*}\mathbb{R}^{n},
\bar{U}\Subset\tilde{U}\Subset T^{*}\mathbb{R}^{n}, \tilde{U}, \tilde{V}
\text{ open }
\end{equation*}
\begin{equation*}
A\left(T-T'\right)B=\mathcal{O}\left(h^{\infty}\right)\colon\mathcal{D}_{h}'
\left(\mathbb{R}^{n}\right)\rightarrow C^{\infty}\left(\mathbb{R}^{n}\right).
\end{equation*}
We shall also use the notation $T\equiv T'.$

Lastly, we define global semi-classical Fourier integral operators.
\begin{Def}\label{dfio}
Let $\Lambda\subset
T^{*}\mathbb{R}^{k}$ be a smooth closed
Lagrangian submanifold with respect to the canonical symplectic
structure on $T^{*}\mathbb{R}^{k}.$
Let $r\in\mathbb{R}.$
Then the space $I^{r}_{h}\left(\mathbb{R}^{k}, \Lambda\right)$ of 
semi-classical
Fourier integral
distributions of order $r$ associated to $\Lambda$ is defined as the set
of all $u\in\mathcal{D}'_{h}\left(\mathbb{R}^{k}\right)$
such
that
\begin{equation*}
\left(\prod_{j=0}^{N}
A_{j}\right)\left(u\right)=\mathcal{O}_{L^{2}\left(\mathbb{R}^{k}\right)}\left(h^{N-r-\frac{k}{4}}\right),
h\to 0,
\end{equation*}
for all $N\in\mathbb{N}_{0}$ and for all $A_{j}\in \Psi_{h}^{0}\left(1, 
\mathbb{R}^{k}\right),$ $j=0, \dots, N-1,$ with
compactly
supported symbols and principal symbols vanishing on $\Lambda$, and any $
A_N \in
\Psi_h^{0} ( 1 , \mathbb{R}^{k}) $ with a compactly supported symbol.

A continuous linear operator
$C_{c}^{\infty}\left(\mathbb{R}^{m}\right)\rightarrow\mathcal{D}_{h}'\left(\mathbb{R}^{l}\right),$
whose Schwartz kernel is an element of
$I_{h}^{r}(\mathbb{R}^{m+l}, \Lambda)$ for some
Lagrangian submanifold $\Lambda\subset T^{*}\mathbb{R}^{m+l}$ and some 
$r\in\mathbb{R}$
will be called a global semi-classical Fourier integral
operator of order $r$ associated to $\Lambda.$
We denote the space of these operators by
$\mathcal{I}_{h}^{r}(\mathbb{R}^{m+l}, \Lambda).$
\end{Def}

\end{document}